\documentclass{article}
\usepackage{amsfonts}
\usepackage{amsmath}
\usepackage{amsmath,amssymb,latexsym,color}
\usepackage[mathscr]{eucal}
\newtheorem{thm}{Theorem}[section]

\newtheorem{lemma}[thm]{Lemma}

\newtheorem{example}{Example}[section]
\newtheorem{defin}[thm]{Definition}

\newtheorem{fact}{Fact}[section]

\newcommand{\proof}{{\it Proof.\quad}}
\newcommand{\qed}{\hfill\Box\medskip}
\usepackage{CJK}
\begin{document}
\begin{CJK*}{GBK}{song}

\renewcommand{\baselinestretch}{1.3}
\title{\bf Full automorphism groups of association schemes based on attenuated spaces}
\author{
Wen Liu$^{\rm a,b}$\quad Kaishun Wang$^{\rm a}$\footnote{Corresponding author. E-mail address: wangks@bnu.edu.cn} \\
{\footnotesize  $^{\rm a}$ \em  Sch. Math. Sci. {\rm \&} Lab. Math.
Com. Sys.,
Beijing Normal University, Beijing, 100875,  China}\\
\footnotesize $^{\rm b}$ \em  Math. {\rm \&} Inf.  Sci.  College,
Hebei Normal University, Shijiazhuang, 050024, China}
\date{}
\maketitle

\begin{abstract}
The set of  subspaces with a given dimension in an attenuated space
has a structure of a symmetric association scheme, which is a
generalization of both Grassmann schemes and bilinear forms schemes.
In [K. Wang, J. Guo, F. Li, Association schemes based on attenuated
space, European J. Combin. 31 (2010) 297--305], its intersection
numbers were computed. In this paper, we determine its full
automorphism group.

\medskip
\noindent 2010 {\em AMS classification:} 05E30

\end{abstract}

\section{Introduction}

Let $\mathbb{F}_q$ be a finite field with $q$ elements, where $q$ is
a prime power. For two non-negative integers $n$ and $l$, suppose
$\mathbb{F}_q^{(n+l)}$   denotes an $(n+l)$-dimensional   row vector
space  over $\mathbb{F}_q$. The set of all  matrices over
$\mathbb{F}_q$ of the form
$$
\left(\begin{array}{cc}
 T_{11} & T_{12}\\
0 & T_{22}
\end{array}\right),$$
where $T_{11}$ and $T_{22}$ are nonsingular $n\times n$ and $l\times
l$ matrices  respectively, forms a group under matrix
multiplication, called the {\it singular general linear group} of
degree $n+l$ over $\mathbb{F}_q$ and denoted by
$GL_{n+l,\,n}(\mathbb{F}_q)$.

Let $P$ be an $m$-dimensional subspace of $\mathbb{F}_q^{(n+l)}$,
denote also by $P$ an $m\times (n+l)$ matrix of rank $m$ whose rows
span the subspace $P$ and call the matrix $P$ a matrix
representation of the subspace $P$. The group
$GL_{n+l,\,n}(\mathbb{F}_q)$ acts on $\mathbb{F}_q^{(n+l)}$ by the
vector matrix multiplication. This action induces an action on the
set of subspaces of $\mathbb{F}_q^{(n+l)}$, i.e., a subspace $P$ is
carried by $T\in GL_{n+l,\,n}(\mathbb{F}_q)$ into the subspace $PT$.
The vector space $\mathbb{F}_q^{(n+l)}$ together with this group
action  is called the $(n+l)$-dimensional {\it singular linear
space} over $\mathbb{F}_q$. This concept was introduced in
\cite{wgl1, wgl2}.

For $1\leq i\leq n+l$, let $e_i $ be the vector in
$\mathbb{F}_q^{(n+l)}$ whose $i$-th coordinate is 1 and all other
coordinates are 0. Denote by $E$ the $l$-dimensional subspace of
$\mathbb{F}_q^{(n+l)}$ generated by
$e_{n+1},e_{n+2},\ldots,e_{n+l}$. An $m$-dimensional subspace $P$ of
$\mathbb{F}_q^{(n+l)}$ is called a
 subspace of {\it type} $(m,k)$ if $\dim(P\cap E)=k$. The collection
of all the subspaces of type $(m,0)$ in $\mathbb{F}_q^{(n+l)}$ is
the {\em attenuated space}, denoted by $A_q(n, l)$.

 Denote by $X_m$
the set  of all subspaces of type $(m,0)$ in $\mathbb{F}_q^{(n+l)}$.
For $0\leq i\leq {\mathrm{min}}\{m,n-m\}$ and $0\leq j-i\leq
{\mathrm{min}}\{m-i,j,l\},$ let
$$ R_{i,j-i}=\{(P,Q)\in X_m\times X_m |\;
{\mathrm{dim}}(P'\cap Q')=m-i,\,{\mathrm{dim}}(P\cap Q)=m-j\},
$$
where
\begin{equation}\label{PQ}
P=\bordermatrix{ &_n&_l\cr &P'&P''},\quad Q=\bordermatrix{&_n&_l\cr
&Q'&Q''}.
\end{equation}
 Wang et al. \cite{wgl1} proved that the configuration
$$ \mathfrak{X}_m=(X_m,
\{R_{i, j-i}\}_{0\leq i\leq {\mathrm{min}}\{m,n-m\}, 0\leq j-i\leq
{\mathrm{min}}\{m-i,j,l\}})
$$
 is a symmetric association scheme,
and computed its intersection numbers. Recently Kurihara
\cite{Kurihara} computed its character table. If $l=0$, the   scheme
$\mathfrak{X}_m$ is the Grassmann scheme $J_q(n,m)$; and if $m=n$,
the   scheme $\mathfrak{X}_m$  is the bilinear forms scheme
$H_q(n,l)$. We refer  readers to  \cite{Bannai, bcn} for the general
theory of association schemes.

Let $\mathfrak{X}=(X, \{R_i\}_{0\leq i\leq d})$ be an association
scheme. If a permutation $\sigma$ on $X$ induces a permutation
$\bar\sigma$ on $\{0, 1,\ldots, d\}$ by $(\sigma(x), \sigma(y))\in
R_{\bar\sigma(i)}$ for $(x, y)\in R_i$, then $\sigma$ is called an
{\em automorphism} of $\mathfrak{X}$. The set of all automorphisms
of $\mathfrak{X}$ becomes a group, called the {\em full automorphism
group} of $\mathfrak{X}$, denoted by ${\rm Aut} (\mathfrak{X}).$ An
automorphism of $\mathfrak{X}$ is called an {\em inner automorphism}
if it induces the identity permutation on $\{0,1,\ldots, d\}$.
Clearly, the set of all inner automorphisms of $\mathfrak{X}$ is a
normal subgraph of ${\rm Aut} (\mathfrak{X}),$ which is called the
{\em inner automorphism group} of $\mathfrak{X}$, denoted by ${\rm
Inn} (\mathfrak{X})$.

In 1949, Chow   determined the full automorphism group  of the
Grassmann scheme $J_q(n,m)$.
\begin{thm}{\rm(\cite{chow})} \label{thm:chow}
Let $1<m<n-1.$ Then
$$
{\rm Aut} (J_q(n,m))= \left\{
                            \begin{array}{ll}
                              P\Gamma L(n,\mathbb {F}_{q}), & \hbox{if $n\neq 2m,$} \\
                              P\Gamma L(n,\mathbb
{F}_{q}).2, & \hbox{if $n=2m$.}
                            \end{array}
                          \right.
$$
\end{thm}

In \cite{koolen2}, Fujisaki et al. determined the full automorphism
group of the twist Grassmann scheme $\tilde J_q(2e+1, e)$. This
scheme has the same parameters as $J_q(2e+1, e)$; see
\cite{koolen1}. In 1965, Deng and Li determined the full
automorphism group of the bilinear forms scheme $H_q(n,l)$.
\begin{thm}{\rm(\cite{deng})} \label{thm:deng}
Let $n$ and $l$ be two integers not less than $2$. Then
$$
{\rm Aut} (H_q(n,l))= \left\{
                            \begin{array}{ll}
                              P\Gamma L(n+l,\mathbb
{F}_{q})_{E}, & \hbox{if $n\neq l$,} \\
                              P\Gamma
L(n+l,\mathbb {F}_{q})_{E}.2, & \hbox{if $n=l$.}
                            \end{array}
                          \right.
$$
\end{thm}

Observe that $\mathfrak{X}_1$ is an association scheme with two
classes. Since one   relation graph  is  ${n\brack 1}$ copies of the
complete graph on ${q^l}$ vertices, it's full automorphism group is
the wreath product $S_{q^l}\wr S_{n\brack 1}$.

Motivated by above results,  in this paper we shall determine the
full automorphism group of  $\mathfrak{X}_m$, and obtain the
following result.

\begin{thm}\label{thm:main}
Let $1<m<n-1$ and $l>0$. Then    ${\rm Aut}(\mathfrak{X}_m)= P\Gamma
L(n+l,\mathbb {F}_{q})_E$.
\end{thm}

\section{Proof of Theorem~\ref{thm:main}}

In this section we always assume that   $1<m<n-1$ and $l>0$. For
each integer $k$ with $2\leq k\leq m$, let $\Gamma^{(k)}$ denote the
relation graph $(X_k, R_{1, 0})$ of $\mathfrak{X}_k$. We first
determine the full automorphism group of $\Gamma^{(m)}$, then prove
Theorem~\ref{thm:main}.

Note that  an $m$-dimensional subspace   $P$ of form (\ref{PQ}) in
$\mathbb F_q^{(n+l)}$ is a vertex of $\Gamma^{(m)}$ if and only if
${\rm rank} (P')=m.$ Therefore,  two vertices $P,Q$ of
$\Gamma^{(m)}$ are adjacent if and only if their sum  $P+Q$  is a
subspace of type $(m+1, 0)$.

\begin{lemma}\label{ll2}
Let $P$ and $Q$  be two vertices as in (\ref{PQ}) of $\Gamma^{(m)}$.
If ${\rm dim}(P'\cap Q')=m-i$ and ${\rm dim}(P\cap Q)=m-j$, then the
distance of $P$ and $Q$
$$\partial(P,Q)
= \left\{\begin{array}{ll} j,& \hbox{if $ i>0,$}\\
j+1,& \hbox{if $i=0.$}\end{array}\right.
$$
\end{lemma}

\proof In the Grassmann graph $J_q(n+l,m)$,  two vertices $x$ and
$y$ are at distance $j$ if and only if $\dim (x\cap y)=m-j.$  Since
$\Gamma^{(m)}$ is a subgraph of  $J_q(n+l,m)$, by ${\rm dim}(P\cap
Q)=m-j$ we have $\partial(P,Q)\geq j$. Write $P\cap
Q=W=\bordermatrix{ &_n&_l\cr &W'&W''}$ and
$$
P=\left(
\begin{array}{cc}
W'&W''\\
   \alpha_{1}' & \alpha_{1}''\\
 \vdots & \vdots\\
 \alpha_{j-i}' & \alpha_{j-i}''\\
\delta_{1}' & \delta_{1}''\\
 \vdots & \vdots\\
  \delta_{i}' & \delta_{i}''
  \end{array}
              \right),
\quad Q=\left(
\begin{array}{cc}
W'&W''\\

\alpha_{1}' & \beta_{1}''\\
 \vdots & \vdots\\
 \alpha_{j-i}' &\beta_{j-i}''\\
 \gamma_{1}' & \gamma_{1}''\\
  \vdots & \vdots\\
   \gamma_{i}' & \gamma_{i}''
  \end{array}
              \right).
$$

{\bf Case 1}\,\ $i>0.$

If $j=i$, write
 $$
 P_{1}=\left(
\begin{array}{cc}
  W' &W''\\
  \gamma_{1}' & \gamma_{1}''\\
  \delta_{2}' & \delta_{2}''\\
 \delta_{3}' & \delta_{3}''\\
  \vdots & \vdots\\
  \delta_{i}' & \delta_{i}''
 \end{array}\right),\;
P_{2}=\left(
\begin{array}{cc}
   W' &W''\\
  \gamma_{1}' & \gamma_{1}''\\
  \gamma_{2}' & \gamma_{2}''\\
  \delta_{3}' & \delta_{3}''\\
  \vdots & \vdots\\
  \delta_{i}' & \delta_{i}''
  \end{array}\right),\;
\ldots, P_{j-1}= \left(
\begin{array}{cc}
   W' &W''\\
   \gamma_{1}' & \gamma_{1}''\\
\gamma_{2}' & \gamma_{2}''\\
  \vdots & \vdots\\
  \gamma_{i-1}' & \gamma_{i-1}''\\
  \delta_{i}' & \delta_{i}''
 \end{array}\right),
$$
then $(P, P_1,\ldots, P_{j-1}, Q)$ is a path of length $j$ from $P$
to $Q$. Therefore $\partial(P,Q)=j$.

If $j>i$, write
$$
P_1=\left(
\begin{array}{cc}
  W' &W''\cr
  \alpha_{1}' & \alpha_{1}''\cr
  \vdots & \vdots\cr
  \alpha_{j-i-2}' & \alpha_{j-i-2}''\cr
   \alpha_{j-i-1}' & \alpha_{j-i-1}''\cr
   \gamma_1' & \gamma_1''\cr
   \delta_{1}' & \delta_{1}''\cr
  \vdots & \vdots\cr
   \delta_{i}' & \delta_{i}''
 \end{array}\right),\;
  P_2=\left(\begin{array}{cc}
   W' &W''\cr
   \alpha_{1}' & \alpha_{1}''\cr
  \vdots & \vdots\cr
   \alpha_{j-i-2}' & \alpha_{j-i-2}''\cr
   \alpha_{j-i}' & \beta_{j-i}''\cr
  \gamma_1' & \gamma_1''\cr
   \delta_{1}' & \delta_{1}''\cr
  \vdots & \vdots\cr
   \delta_{i}' & \delta_{i}''
\end{array}\right),\ldots,
 P_{j-i}=\left(\begin{array}{cc}
   W' &W''\cr
   \alpha_{2}' & \beta_{2}''\cr
  \vdots & \vdots\cr
   \alpha_{j-i-1}' & \beta_{j-i-1}''\cr
   \alpha_{j-i}' & \beta_{j-i}''\cr
   \gamma_1' & \gamma_1''\cr
   \delta_{1}' & \delta_{1}''\cr
  \vdots & \vdots\cr
   \delta_{i}' & \delta_{i}''
  \end{array}\right),
$$
  then $(P, P_1,\ldots, P_{j-i})$ is a path of length $j-i.$ Since
  $  {\rm dim}(P_{j-i}'\cap Q')={\rm dim}(P_{j-i}\cap Q)=m-i,$
  we have $\partial(P_{j-i},Q)=i$. It follows that
 $\partial(P,Q)=j$.

\medskip
{\bf Case 2} \,\ $i=0$.

The neighborhood of $P$ consists of vertices
$$
R=\left(\begin{array}{cc}
\bar{P}'&\bar{P}''\cr
\xi'&\xi''
\end{array}\right),
$$
where
 $\bar{P}=\left(\begin{array}{cc}\bar{P}'& \bar{P}''\end{array}\right)$ is an $(m-1)$-dimensional
 subspace of $P$ and $\xi'\in \mathbb F^{(n)}\setminus  P'$.

If $W\subseteq \bar{P}$, then ${\rm dim}(R\cap Q)=m-j$ and $ {\rm
dim}\bigg(\left(
       \begin{array}{c}
          \bar{P}' \\
         \xi'\\
       \end{array}
     \right)
\cap Q'\bigg)=m-1$. By Case 1, one gets $\partial(R,Q)=j$. If $W
\nsubseteq\bar{P}$, similarly  we have $\partial(R,Q)=j+1$. Hence,
$\partial(P,Q)=j+1$. $\qed$

Next we shall study the first and   second subconstituents of
$\Gamma^{(m)}$. Lemma 2.1 in \cite {wgl1} implies that $P\Gamma L
(n+l, \mathbb F_q)_E$ acts transitively on each $R_{i, j}$, in
particular
 the graph $\Gamma^{(m)}$ is   vertex-transitive. So
  we only need to consider the subconstituents with respect to the
vertex
$$ M =\bordermatrix{&m&n-m&l\cr &I&0&0} =\bordermatrix{&n&l\cr
 &M'&M''}.
$$

Let $ \Gamma_{i,j-i}(M)$ be the set of vertices $P$ of
$\Gamma^{(m)}$ satisfying  $(M, P)\in R_{i, j-i}$. By Lemma
\ref{ll2}  the second subconstituent $\Gamma_{2}(M)$ has a partition
$$
\Gamma_{0,1}(M)\;\dot{\cup}\;\Gamma_{1,1}(M)\;\dot{\cup}\; \Gamma_{2,0}(M).
$$

For simplicity, write
$$
   P(U;\alpha,\beta,\gamma)=\bordermatrix{&m&n-m&l\cr
   &U& 0 & 0  \cr
  &\alpha&\beta & \gamma
   }\begin{array}{c}
   t\cr
   1\end{array},
$$
$$
 P(U;\alpha_1,\beta_1,\gamma_1;\alpha_2,\beta_2,\gamma_2)=\bordermatrix{&m&n-m&l\cr
   &U& 0 & 0  \cr
  &\alpha_1&\beta_1 & \gamma_1\cr
    &\alpha_2&\beta_2& \gamma_2
   }\begin{array}{c}
   t\cr
   1\cr
   1\end{array},
$$
where ${\rm rank} (U)=t$, ${\rm rank} (\begin{array}{cc}\beta &
\gamma\end{array})=1$ and ${\rm rank} \left(\begin{array}{cc}\beta_1
& \gamma_1\cr \beta_2 & \gamma_2
\end{array}\right)=2$.
 Then

$$
\begin{array}{rcl}
\Gamma_{1}(M)&=&\{P(W;\alpha,\beta,\gamma)\mid {\rm rank}(W)=m-1,
\beta\neq 0\},\\
$$\Gamma_{0,1}(M)&=&\{P(W;\alpha,0,\gamma)\mid {\rm rank}(W)=m-1,\alpha\notin W, \gamma\neq
0\},\\
$$\Gamma_{1,1}(M)&=&\{P(U;\alpha_1,\beta_1,\gamma_1;\alpha_2,0,\gamma_2)\mid {\rm rank}(U)
=m-2,\alpha_2\notin U, \beta_1\neq 0,\gamma_2\neq 0\},\\
$$\Gamma_{2,0}(M)&=&\{P(U;\alpha_1,\beta_1,\gamma_1;\alpha_2,\beta_2,\gamma_2)\mid
 {\rm rank}(U)=m-2,{\rm rank}\left({\beta_1\atop \beta_2}\right)=2\}.
\end{array}$$

\begin{lemma}\label{lemma:wan} {\rm (\cite[Corollary 1.9]{wan1})}
Let $0\leq k\leq m\leq n.$ Then the number of $m$-dimensional
subspaces containing a given $k$-dimensional subspace in $\mathbb
F^{(n)}$ is equal to ${n-k\brack m-k}$.
\end{lemma}

\begin{lemma}\label{114}
 {\rm(i)} Let $P(W;\alpha,0,\gamma)\in \Gamma_{0,1}(M)$ and
$P(W';\alpha',\beta',\gamma')\in \Gamma_1(M)$.  Then the two
vertices are adjacent in $\Gamma^{(m)}$ if and only if $W=W'.$ In
particular, each vertex in $\Gamma_{0,1}(M)$ has $q^{l+1}{n-m\brack
1}$ neighbors  in $\Gamma_1(M)$;

 {\rm(ii)} Each vertex in $\Gamma_{1,1}(M)$ has  $q^2$ neighbors  in $\Gamma_1(M)$;

 {\rm(iii)} Each vertex in $\Gamma_{2,0}(M)$ has  $(q+1)^2$ neighbors  in
$\Gamma_1(M)$.
\end{lemma}
\proof (i)   Suppose $P(W;\alpha,0,\gamma)$ is adjacent to
$P(W';\alpha',\beta',\gamma')$. If $W\neq W'$, by $\beta'\neq 0$ and
$\gamma\neq 0$ the dimension of
$P(W;\alpha,0,\gamma)+P(W';\alpha',\beta',\gamma')$ is  $m+2$,  a
contradiction. The converse is immediate from $\alpha\notin W$ and
$\beta'\neq 0$.  Therefore, the first statement is valid.

 Since $P(W;\alpha',\beta',\gamma')$  is of
type $(m, 0)$, the rank of $\left(
      \begin{array}{cc}
        W & 0 \\
        \alpha' & \beta'\\
      \end{array}
    \right)
$ is $m$. By Lemma~\ref{lemma:wan} the subspace  $\left(
      \begin{array}{cc}
        W & 0 \\
        \alpha' & \beta'\\
      \end{array}
    \right)
$ with $\beta'\neq 0$ has $q{n-m\brack 1}$ choices. For a given
subspace  $\left(
      \begin{array}{cc}
        W & 0 \\
        \alpha' & \beta'\\
      \end{array}
    \right),
$  there are $q^l$ choices for $\gamma'$. Hence,  (i) holds.

(ii) Given a vertex
$P(U;\alpha_1,\beta_1,\gamma_1;\alpha_2,0,\gamma_2)\in
\Gamma_{1,1}(M)$. We claim that this vertex is adjacent to
$P(W;\alpha,\beta,\gamma)\in \Gamma_1(M)$ if and only if $U\subseteq
W $ and $P(W;\alpha_1,\beta_1,\gamma_1;\alpha_2,0,\gamma_2)$ is a
subspace of type $(m+1,0)$ containing $P(W;\alpha,\beta,\gamma)$.
Suppose $P(W;\alpha,\beta,\gamma)$ is adjacent to
$P(U;\alpha_1,\beta_1,\gamma_1;\alpha_2,0,\gamma_2)$. If
$U\nsubseteq W$, by $\beta_1\neq 0$ and $\gamma_2\neq 0, $  the
dimension of
$P(W;\alpha,\beta,\gamma)+P(U;\alpha_1,\beta_1,\gamma_1;\alpha_2,0,\gamma_2)$
is at least $m+2$, a contradiction. It follows  that $U\subseteq W$.
Then
$$
P(W;\alpha,\beta,\gamma)+P(U;\alpha_1,\beta_1,\gamma_1;\alpha_2,0,\gamma_2)=
P(W;\alpha,\beta,\gamma)+P(W;\alpha_1,\beta_1,\gamma_1;\alpha_2,0,\gamma_2),
$$
and so $P(W;\alpha_1,\beta_1,\gamma_1;\alpha_2,0,\gamma_2)$ is a
subspace of type $(m+1,0)$ containing $P(W;\alpha,\beta,\gamma)$.
The converse is immediate. Hence, our claim is valid.

  By
Lemma~\ref{lemma:wan} the number of $(m-1)$-subspaces $W$ in
$\mathbb F^{(m)}$ containing $U$ is $q+1$. Observe that
$P(W;\alpha_1,\beta_1,\gamma_1;\alpha_2,0,\gamma_2)$ is of type
$(m+1, 0)$ if and only if $W\neq\left(\begin{array}{c} U\cr \alpha_2
\end{array}
\right)$. Then $W$ has $q$ choices. For a fixed $W$, by
Lemma~\ref{lemma:wan} again  $P(W;\alpha,\beta,\gamma)$ has $q$
choices. Therefore, (ii) holds.

The proof of (iii)  is similar to that of (ii), and will be omitted.
$\qed$

For a   subspace $W$ of type $(m-1,0)$ in $\mathbb{F}_q^{(n+l)}$,
let $C(W)$ be the set of all vertices of
 $\Gamma^{(m)}$ containing  $W$. For convenience we  denote by $\Delta$
the induced subgraph of $\Delta$ on $\Gamma^{(m)}$.

\begin{lemma}\label{115}
Let $W$ be a subspace of type $(m-1,0)$ in $\mathbb{F}_q^{(n+l)}$.
Then  $C(W)$ is isomorphic to the complete multipartite graph
$K_{n-m+1\brack 1}(q^l)$.
\end{lemma}

\proof Let $W=\bordermatrix{ &_n&_l\cr &W'&W''}$. Then  $C(W)$
consists  of vertices
\begin{eqnarray}\label{a0}
\left(\begin{array}{cc}
W'&W''\cr
\alpha'&\alpha''
\end{array}
\right),
\end{eqnarray}
where $ {\rm rank}\left(\begin{array}{c} W'\cr \alpha'
\end{array}
\right)=m.$ By Lemma~\ref{lemma:wan}   the number of
$m$-dimensional subspaces $ \left(\begin{array}{c} W'\cr \alpha'
\end{array}
\right) $ in  $\mathbb{F}_q^{(n)}$  is ${n-m+1\brack 1}$. Since
$\alpha''$ has $q^l$ choices, the subgraph $C(W)$ has
$q^l{n-m+1\brack 1}$ vertices.

For a  given $m$-dimensional subspace $ \left(\begin{array}{c} W'\cr
\alpha'
\end{array}
\right) $ in $\mathbb{F}_q^{(n)}$, the vertices of form (\ref{a0})
form an independent set with $q^l$ vertices. Note that all these
independent sets form a partition of  $C(W)$.
 Since each vertex in an
independent set is adjacent to any vertex in the remaining
independent sets, the desired result follows.$\qed$

\begin{lemma}\label{ll1}
 If an induced subgraph $\Delta$ of $\Gamma^{(m)}$ is isomorphic to $K_{n-m+1\brack 1}(q^l)$, then
 $\Delta$ is a subgraph $C(W)$, where $W$ is a subspace of type $(m-1,0)$ in $\mathbb{F}_q^{(n+l)}$.
\end{lemma}
\proof Since $\Gamma^{(m)}$ is vertex-transitive, we may assume that
$\Delta$   contains $M$. Pick $X\in \Delta\cap \Gamma_1(M)$, and
write $M\cap X=W$. Then
$$W=\bordermatrix{&m&n-m&l\cr &W'&0&0},\quad
 X=P(W'; \alpha', \beta', \gamma').
$$

 Now we shall show that $\Delta=C(W)$.
Suppose $\Delta\cap (\Gamma_{1,1}(M)\cup\Gamma_{2,0}(M)) \neq
\emptyset$. Pick a vertex $P$ from this set. Then the vertices $P$
and $M$ have $q^{l+1}{n-m\brack 1}$ common neighbors in $\Delta$, a
contradiction to Lemma~\ref{114} (ii), (iii).  So  $\Delta\cap
\Gamma_{0,1}(M)\neq \emptyset$. Pick a vertex
$P(W'';\alpha'',0,\gamma'')\in \Delta\cap\Gamma_{0,1}(M)$. Since
$\Delta$ is a  complete multipartite graph, the vertices $X$ and
$P(W'';\alpha'',0,\gamma'')$ are adjacent. By
 Lemma \ref{114} (i), one gets $W'=W''$. It
follows that $\Delta\cap\Gamma_{0,1}(M)\subseteq C(W)$. Similarly,
we have $\Delta\cap\Gamma_1(M)\subseteq C(W)$. Hence,
$\Delta\subseteq C(W)$. By Lemma \ref{115} the subgraphs $\Delta$
and $C(W)$ have the same number of vertices, so $\Delta=C(W)$.
$\qed$

\begin{lemma}\label{ll6}
${\rm Aut}(\Gamma^{(m)})=P\Gamma L(n+l,\mathbb {F}_q)_E$.
\end{lemma}
\proof We first prove that the result holds for $m=2.$ Pick any
automorphism $\tau$ of $\Gamma^{(2)}$. Then  $\tau$ is a permutation
 on the set of lines  and permutes the points of the   attenuated space $A_q(n, l)$. By Deng and
Li's result in \cite{deng}, the automorphism $\tau$ can be extended
to a collineation fixing $E$ of the projective space $PG(n+l,
\mathbb F_q)$. By  the fundamental theorem of the projective
geometry \cite[Theorem 2.23]{wan3}, we have $\tau\in P\Gamma
L(n+l,\mathbb {F}_q)_E.$
 Thus ${\rm Aut}(\Gamma^{(2)})\subseteq
P\Gamma L(n+l,\mathbb {F}_q)_E$. On the other hand, $P\Gamma
L(n+l,\mathbb {F}_q)_E\subseteq {\rm Aut}(\Gamma^{(2)})$. Hence,
${\rm Aut}(\Gamma^{(2)})= P\Gamma L(n+l,\mathbb {F}_q)_E$.

Now let $m\geq 3 $ and $\sigma_m$ be an automorphism of the graph
$\Gamma^{(m)}$. By Lemmas \ref{115} and   \ref{ll1}, the
automorphism $\sigma_m$ induces a permutation  on the set
$\{C(W)\mid W\in X_{m-1}\}$, and  further induces  a permutation
$\sigma_{m-1}$ on $X_{m-1}$. For any two adjacent vertices
$W_{m-1},W_{m-1}'$ of $\Gamma^{(m-1)}$, we have $C(W_{m-1})\cap
C(W_{m-1}')\neq \emptyset$. Therefore $\sigma_m(C(W_{m-1})\cap
C(W_{m-1}'))=C(\sigma_{m-1}(W_{m-1}))\cap
C(\sigma_{m-1}(W_{m-1}'))\neq \emptyset$, which implies that the two
vertices $\sigma_{m-1}(W_{m-1})$ and $\sigma_{m-1}(W_{m-1}')$ are
adjacent in $\Gamma^{(m-1)}$. Hence, $\sigma_{m-1}\in {\rm Aut}
(\Gamma^{(m-1)})$.

 By induction, for each $3\leq k\leq m,$ the map
$f_k: \sigma_k\longmapsto \sigma_{k-1}$ is a  homomorphism from
${\rm Aut}(\Gamma^{(k)})$ to ${\rm Aut}(\Gamma^{(k-1)})$. We claim
that $f_k$ is injective. Suppose $\sigma_{k-1}=\iota$, the identity
permutation on $X_{k-1}$. For each vertex $W_{k}$ of $\Gamma^{(k)}$,
there exist two vertices $W_{k-1}$ and $ W_{k-1}'$ of
$\Gamma^{(k-1)}$ such that $W_{k}=W_{k-1}+ W_{k-1}'$. Since
$\{W_{k}\}=C(W_{k-1})\cap C(W_{k-1}'),$ we have $\{\sigma_k
(W_{k})\}=C(\sigma_{k-1}(W_{k-1}))\cap
C(\sigma_{k-1}(W_{k-1}'))=\{W_{k}\}$. Hence our claim is valid. It
follows that   $|{\rm Aut}(\Gamma^{(m)})|\leq|{\rm
Aut}(\Gamma^{(2)})|.$ Since $ P\Gamma L(n+l,\mathbb {F}_q)_E$ is a
subgroup of ${\rm Aut}(\Gamma^{(m)})$, the desired result follows.
$\qed$

\medskip{\em Proof of Theorem~\ref{thm:main}:} \,\ The fact that $P\Gamma L (n+l, \mathbb
F_q)_E$ acts transitively on each $R_{i, j-i}$ implies   ${\rm
Inn}(\mathfrak{X}_m)={\rm Aut}(\Gamma^{(m)}).$ In order to show our
result, it suffices to prove that all  valencies $n_{i, j-i}$'s of
$\mathfrak{X}_m$ are pairwise distinct.

 By \cite{wgl1} we have
$$
n_{i,j-i}=q^{i^2+il+\frac{(j-i)(j-i-1)}{2}}{n-m\brack i} {m \brack
i}{m-i\brack j-i} \prod \limits_{s=l-(j-i)+1}^l(q^s-1).
$$
Suppose $n_{a, b-a}= n_{a', b'-a'}$. Since $q$ is a prime power, we
obtain
\begin{eqnarray}\label{a1}
a^2+al+\frac{(b-a)(b-a-1)}{2}=a'^2+a'l+\frac{(b'-a')(b'-a'-1)}{2},
\end{eqnarray}
\begin{eqnarray}\label{a2}
{n-m\brack a}{m \brack a}{m-a\brack b-a} \prod
\limits_{s=l-(b-a)+1}^l(q^s-1)={n-m\brack a'} {m \brack
a'}{m-a'\brack b'-a'} \prod \limits_{s=l-(b'-a')+1}^l(q^s-1).
\end{eqnarray}
Simplifying (\ref{a1}), we have
\begin{eqnarray}\label{a3}
2(a-a')(a+a'+l)=((b'-a')-(b-a))(b'-a'+b-a-1).
\end{eqnarray}

Let
$$\begin{array} {rcl}
f_{i,j}(x)&=&\prod\limits_{s=1}^{i}(x^s-1)^2\prod\limits_{s=1}^{j-i}(x^s-1).\\
g_{i,j}(x)&=&\prod\limits_{s=n-m-i+1}^{n-m}(x^s-1)\prod\limits_{s=m-i+1}^{m}(x^s-1)\prod\limits_{s=m-j+1}^{m-i}(x^s-1)
\prod\limits_{s=l-j+i+1}^{l}(x^s-1),
\end{array}
$$
The equality (\ref{a2}) implies that
$f_{a,b}(q)g_{a',b'}(q)=f_{a',b'}(q)g_{a,b}(q)$  for all prime
powers $q$, so $f_{a,b}(x)g_{a',b'}(x)=f_{a',b'}(x)g_{a,b}(x).$
Since $1$ is a root with multiplicity $2a+(b-a)+2a'+2(b'-a')$ of
$f_{a,b}(x)g_{a',b'}(x)$ and $1$ is a root with multiplicity
$2a+2(b-a)+2a'+(b'-a')$ of $f_{a',b'}(x)g_{a,b}(x)$, we obtain
 $b'-a'=b-a$. By (\ref{a3}) one gets $a'=a$,  and so  $b'=b$. Hence
the desired result follows.$\qed$

\section*{Acknowledgement}
W. Liu's research is supported by NSFC(11171089, 11271004). K.
Wang's research is supported by  NSFC(11271047), SRFDP and the
Fundamental Research Funds for the Central University of China.

\end{CJK*}

\end{document}